\documentstyle[12pt,dina4]{article}
 
\newtheorem{THEOREM}{Theorem}[section]
\newcommand{\thm}[1]{
\begin{THEOREM}
#1
\end{THEOREM}
}
\newtheorem{rem}{Remark}[section]

\newtheorem{LEMMA}[THEOREM]{Lemma}
\newcommand{\la}[1]{
\begin{LEMMA}
#1
\end{LEMMA}
} 
\newtheorem{COROLLARY}[THEOREM]{Corollary}
\newcommand{\coro}[1]{
\begin{COROLLARY}
#1
\end{COROLLARY}
}
\newtheorem{PROPOSITION}[THEOREM]{Proposition}

\newtheorem{DEFINITION}{Definition}[section]
 
\newcommand{\nm}[1]{\parallel #1 \parallel}

\newcommand{\nmwg}[1]{\nm{#1}_{L^p_{v(\gamma )}}}

\begin{document}
 \begin{center}
 {\LARGE{}Fractional integration for Laguerre expansions}\\[.4cm]
 {\sc{}George Gasper \footnote{Dept. of Mathematics, Northwestern University, 
 Evanston, IL 60208, USA. Supported in part by NSF Grant DMS--9103177.},
 Krzysztof Stempak \footnote{Institute of Mathematics, 
 University of Wroc\l aw, Pl. Grunwaldzki 2/4, 50 --384 
Wroc\l aw, Poland. Research of this author was 
 done during a stay at the Fb Mathematik of the TH Darmstadt and was 
 supported by the Deutsche Forschungsgemeinschaft under grant 436POL/115/1/0.}
 and Walter  Trebels \footnote{Fb. Mathematik, TH Darmstadt, 
 Schlo\ss{}gartenstr. 7, D--64289 Darmstadt, Germany.} }\\[.4cm]
  {(July 8, 1994 version)}\\[.4cm]
 \end{center}


\bigskip 
{\bf Key words.} Fractional integration, multiplier conditions, Laguerre 
polynomials 

\bigskip
{\bf AMS(MOS) subject classifications.} 33C45, 42A45, 42C10

\bigskip 
\section{ Introduction}
The aim of this note is to provide a fractional integration theorem 
in the framework of Laguerre expansions. The method of proof consists of 
establishing an asymptotic estimate for the involved kernel and then  
applying a method of Hedberg \cite{pro}. We combine this
result with sufficient $(p,p)$ multiplier criteria of Stempak and Trebels 
\cite{ST}. The resulting sufficient $(p,q)$ multiplier criteria are comparable 
with necessary ones of Gasper and Trebels \cite{laguerre}.

\medskip \noindent 
Our notation is essentially that in \cite{ST}. Thus we consider the Lebesgue spaces
$$L^p_{v(\gamma )} = \{ f: \; \nmwg{f} = ( \int _0^\infty
|f(x) |^p x^{\gamma }\, dx )^{1/p} < \infty \} \; ,\quad 1 \le p <
\infty  ,\quad \gamma >-1,$$
and define the  Laguerre function system $\{ l_k^\alpha \} $ by
$$l_k^\alpha (x) =( k!/\Gamma (k+\alpha +1) ) ^{1/2} e^{-x/2}L_k^\alpha (x),
\quad \alpha >-1,\quad k\in {\bf N}_0. $$
This system is an orthonormal basis in $L^2({\bf R}_+,\, x^\alpha dx)$ and for  
$  \gamma < p(\alpha +1)-1$ we can associate 
to  any $f\in L^p_{v(\gamma )}$  the Laguerre series
$$ f(x) \sim \sum _{k=0}^\infty a_k l_k^\alpha (x),\quad \quad 
a_k = \int _0^\infty f(x)l_k^\alpha (x)x^\alpha  dx. $$ 
It is convenient to introduce the vector space 
$$E=\{ f(x)=p(x)e^{-x/2} \, :\; 0 \le x < \infty ,\; \; p(x) \; {\rm a \; 
polynomial} \} $$
which is dense in $L^p_{v(\gamma )}$. We note that $f\in E$ has only finitely 
many non-zero Fourier-Laguerre coefficients.
Analogous to the definition of the Hardy and Littlewood fractional integral 
operator for Fourier series  (see \cite[Chap. XII, Sec. 8]{zyg}, we define a
fractional integral operator $I_\sigma ,\; \sigma >0,$ for Laguerre expansions 
by 
$$I_\sigma f(x) = \sum _{k=0}^\infty (k+1)^{-\sigma }a_k l_k^\alpha (x) \, , 
\quad \quad    f\in E.$$
Observing that the $l_k^\alpha $ are eigenfunctions with eigenvalues $\lambda 
_k$ of the differential operator
$$L =-\Big( x\frac{d^2}{dx^2} + (\alpha +1)\frac{d}{dx} -\frac{x}{4} \Big)\, , 
\quad \lambda _k = k+ (\alpha +1)/2,$$
(see \cite[(5.1.2)]{szego}) one realizes that $I_1$ is an integral operator 
essentially inverse to $L$ (see the following Remark 2).
As we will see in Section 2 the fractional integral $I_\sigma f$ can be 
interpreted as a twisted convolution of $f$ with a function 
\begin{equation}\label{kernel}
g_\sigma (x) \sim \Gamma (\alpha +1) \sum _{k=0}^\infty (k+1)^{-\sigma } 
L_k^\alpha (x)e^{-x/2} \, .
\end{equation}
 From Theorem 3.1 in \cite[II]{laguerre} it easily follows that 
$g_\sigma \in L^1_{v(\gamma )}$ when $\alpha -\gamma < \sigma $, 
so that, by the convolution theorem of G\"orlich and Markett \cite{indag}, 
$I_\sigma $ extends to a bounded 
operator from $L^p_{v(\gamma )}$ to $L^p_{v(\gamma )}$ when $0\le \alpha p/2 
\le  \gamma \le \alpha ,\; 1\le p \le \infty .$ 
Our main result is

\bigskip \noindent 
\thm{
Let  $\alpha \ge 0,\; 1<p\le q<\infty .$ Assume further that $ 0< \sigma < 
\alpha +1$, $a < (\alpha +1)/p',\; 
b<(\alpha +1)/q,\; a+b \ge 0. $ Then
$$||I_\sigma f ||_{L^q_{v(\alpha -bq)}}\leq C||f ||_{L^p_{v(\alpha +ap)}},
\quad \frac1q=\frac1p-\frac{\sigma -a-b}{\alpha+1}. $$
}

\bigskip \noindent
{\bf Remarks.} 1) Observe that the upper bounds for $a$ and $b$ imply  
lower bounds for $b$ and $a$, respectively.

\smallskip \noindent 
2) By the above it is clear that the sequence 
$\{ (k+1)^{-\sigma }\} $ which generates the fractional integral can be 
replaced by any sequence $\{ \Omega _\sigma (k)\} $
satisfying
$$ \Omega _\sigma (k) = \sum _{j=0}^J c_j\, (k+1)^{-\sigma -j} + 
O((k+1)^{-\sigma -J}) $$ 
for sufficiently large $J$, say $J\ge \alpha +2$, thus in particular obtaining 
the same result for the sequence 
$\{ \Gamma (k+1) /\Gamma (\sigma + k +1)  \} $ -- see \cite{mume}.

\smallskip \noindent
3) A weaker version of Theorem 1.1 (e.g. in the case $a=b=0$) 
can easily be deduced by the following 
argument. By a slight modification of the proof of
Theorem 3.1 in  \cite[II]{laguerre} we have:

\smallskip \noindent
{\it Let $\alpha >-1$ and $N\in {\bf N}_0,\; N> (2\alpha +2)(1/r -1/2) -
1/3$. If $\{ f_k\} $ is a bounded sequence with 
$\lim _{k\to \infty } f_k =0$ and 
$$\sum _{k=0}^\infty (k+1)^{N+(\alpha +1)/r' 
} |\Delta ^{N+1} f_k| \le K_r(f) < \infty ,\quad 1\le r< \infty ,$$
then there exists a function $f\in L^r_{v(\alpha )} $ with
$$ \| f\| _{L^r_{v(\alpha )} } \le C K_r(f),\quad 
f(x) \sim \sum _{k=0}^\infty f_k L_k^\alpha (x)e^{-x/2}. $$
}

\smallskip \noindent
This applied to the sequence $\{ (k+1)^{-\sigma } \} $ gives, by Young's 
inequality (see \cite{indag}),
$$||I_\sigma f ||_{L^q_{v(\alpha )}}\leq C||f ||_{L^p_{v(\alpha )}},
\quad \frac1q>\frac1p-\frac{\sigma }{\alpha+1}, $$
where $\alpha \ge 0,\; \sigma >0$ and $1\le p,\, q \le  \infty $.

\bigskip \noindent
Next we indicate how Theorem 1.1 can be used to gain some insight into 
the structure of $M^{p,q}$--Laguerre multipliers. For the sake of simplicity 
let us restrict ourselves to the case $\gamma = \alpha $. Consider a 
sequence $m = \{ m_k\} $ 
of  numbers and associate to $m$ the operator 
$$ T_m f(x) = \sum _{k=0}^\infty m_ka_k l_k^\alpha (x)\, ,\quad f\in E. $$
The sequence $m$ is called 
a bounded $(p,q)$-multiplier,  notation $ m \in M^{p,q}_{\alpha ,\alpha }$, if
$$   ||\,m\,||_ {M^{p,q}_{\alpha ,\alpha }} := \inf \{ C:\,  ||T_m f||_
{L^q_{v(\alpha ) }}\le C \| f \| _{L^p_{v(\alpha ) }} \quad {\rm for \; all}
\;  f\in E \} $$
is finite. 
If $p \le 2 \le q$, then sufficient conditions follow at once in the following 
way: observe that $M^{2,2}_{\alpha ,\alpha } = l^\infty $, choose 
$\sigma _0 , \, \sigma _1 \ge 0$ such that $I_{\sigma _0}:L^p \to L^2,\; 
I_{\sigma _1}:L^2 \to L^q,\; \sigma _0 +\sigma _1 =\sigma $,
and hence 
$$ ||\,m\,||_ {M^{p,q}_{\alpha ,\alpha }} \le \| k^{-\sigma _0 }\| 
_ {M^{p,2}_{\alpha ,\alpha }} || \{k^\sigma m_k \} ||_ 
{M^{2,2}_{\alpha ,\alpha }} \| k^{-\sigma _1 }\| 
_ {M^{2,q}_{\alpha ,\alpha }}\, ;$$
in particular, $m \in M^{p,q}_{\alpha ,\alpha } $  when $\{k^\sigma m_k \} 
\in l^\infty $. Thus  when $1< p\le 2 \le q < \infty $, and $ \sigma $ is as 
in Theorem 1.1, $\{ (-1)^k (k+1)^{-\sigma } \} \in M^{p,q}_{\alpha ,\alpha } $ 
generates a bounded 
operator which does not fall under the scope of the above introduced fractional 
integral operators. To formulate a Corollary based on a combination of Theorem 
1.1 and the multiplier result in \cite{ST} we need the notion of a 
difference operator $\Delta ^s$ of fractional order $s$ given by
$$ \Delta ^s m_k = \sum _{j=0}^\infty A_j^{-s-1}m_{k+j}\, ,\quad
A_j^t =\frac{\Gamma (j+t+1)}{j!\, \Gamma (t+1)} ,\quad t\in {\bf R}, $$
whenever the sum converges. In view of the remark concerning the case $p<2<q$ 
and on account of duality we may restrict ourselves to the case $ 1<p <q <2$.

\bigskip \noindent
\coro{
If $\alpha \ge 0,\; 1<p <q <2,$ and $s> {\rm max} \{ (2\alpha +2)(1/q -1/2),\, 
1\} $  then, for some constant $C$ independent of 
the sequence $\{ m_k \} $, there holds
$$ \|  \{ m_k \} \| ^2_{M^{p,q}_{\alpha ,\alpha }} \le C \left( 
\|  \{k^\sigma  m_k \} \| ^2_\infty + \sup _N \sum _{k=N}^{2N}
|(k+1)^{s+\sigma } \Delta ^s m_k|^2 (k+1)^{-1} \right) \, .$$
}

\medskip \noindent
The proof follows as in \cite{gtpq}. The result itself should be compared with 
the corresponding necessary condition in \cite[I]{laguerre} which also shows 
that $\{ (-1)^k (k+1)^{-\sigma } \} 
\notin M^{p,q}_{\alpha ,\alpha } $ provided $1<p<q<2$ and $\alpha $ is 
sufficiently large. 

\bigskip \noindent
The plan of the paper is as follows. In Section 2 we derive an asymptotic 
estimate of the  function $g_\sigma $ defined above by (\ref{kernel}). 
Then the twisted generalized convolution is used to dominate $I_\sigma f$ by a 
generalized Euclidean convolution of $g_\sigma $ with $f$. The latter's mapping 
behavior is discussed by a method of Hedberg \cite{pro} which uses maximal 
functions; thus giving Theorem 1.1 in the standard weight case $a=b=0$.
In Section 3 we extend this result to some power weights, modifying an argument 
in Stein and Weiss \cite{sw}. 

\section{Proof of the standard weight case.}
We start by deriving the required asymptotic estimate and showing 
\la{
Let $\alpha >-1$. Then, for $x>0$ and $0<\sigma <\alpha +1$, there holds
\begin{equation}\label{asymptotic}
|g_\sigma (x)|\le C x^{\sigma -\alpha -1}.
\end{equation}
}
{\bf Proof.} 
First note that by subordination  (see \cite[I]{laguerre}, p. 1234)
there holds for $N> \alpha +2,\; N\in {\bf N}, $ 
$$ g_\sigma (x) =C \sum _{k=0}^\infty \left( \Delta ^N(k+1)^{-\sigma } 
\right) L_k^{\alpha +N}(x)e^{-x/2}.$$
Then the assertion (\ref{asymptotic}) follows after it is proved that 
\begin{equation}\label{maj}
\sum _{k=0}^\infty (k+1)^{-\sigma -N}\, |x^{\alpha +1-\sigma } L_k^{\alpha 
+N}(x)e^{-x/2} | \le C
\end{equation}
uniformly in $x>0$. With the notation 
$$  {\cal L}_k^\alpha (x) =( k!/\Gamma (k+\alpha +1)) ^{1/2} 
x^{\alpha /2}e^{-x/2}L_k^\alpha (x) ,\quad k\in {\bf N}_0 .$$
this is equivalent to 
\begin{equation}\label{locmaj}
 \sum _{k=0}^\infty (k+1)^{(\alpha  -N)/2-\sigma 
}\, | {\cal L}_k^{\alpha +N} (x)| \le C x^{\sigma +(N- \alpha -2)/2} 
\end{equation}
uniformly in $0<x<\infty $. We will make use of 
the pointwise estimates for the Laguerre functions in \cite[Sec. 2]{aswa}:
\begin{equation}\label{muest}
|{\cal L}_k^{\alpha + N}(x)| \le C 
\left\{ \begin{array}{l@{\quad if \quad}l}
(x(k+1))^{(\alpha +N)/2} &  0 \le x\le c/(k+1) ,\\
(x(k+1))^{-1/4} &  c/(k+1) \le x \le d(k+1) 
\end{array} \right. 
\end{equation}
for fixed positive constants $c$ and $d$.
These and two further estimates in \cite[(2.5)]{mutrans}    imply that 
\begin{equation}\label{munorm}
\sup _{x>0} |{\cal L}_k^{\alpha + N}(x)| \le C  ,\quad k\in {\bf N}_0,
\end{equation}
and so it is obvious that (\ref{munorm}) implies (\ref{locmaj}) for $x\ge 1 .$ 

\medskip \noindent
Therefore, decomposing dyadically the interval $(0,1)$ it suffices to check 
that 
$$ \sum _{k=0}^\infty (k+1)^{(\alpha  -N)/2-\sigma }\, 
| {\cal L}_k^{\alpha +N} (x)| \le C (2^j)^{\sigma +(N- \alpha -2)/2} $$
provided $2^j \le x \le 2^{j+1},\; j<0$. Using the 
first line of (\ref{muest}) we 
get 
$$ \sum _{k=0}^{2^{-j}}(k+1)^{(\alpha  -N)/2-\sigma }\, 
| {\cal L}_k^{\alpha +N} (x)| \le C (2^j)^{(\alpha +N)/2} 
\sum _{k=0}^{2^{-j}}(k+1)^{\alpha -\sigma }
\le C (2^j)^{\sigma +(N-\alpha -2)/2} ,$$
while the second line of (\ref{muest}) gives
$$ \sum _{k=2^{-j}}^\infty k^{(\alpha  -N)/2-\sigma }\, 
| {\cal L}_k^{\alpha +N} (x)| \le C 2^{-j/4} 
\sum _{k=2^{-j}}^\infty k^{(\alpha  -N)/2-\sigma -1/4}
\le C (2^j)^{\sigma +(N-\alpha -2)/2} .$$
This completes the proof of Lemma 2.1.

\bigskip \noindent
As already mentioned we will apply Hedberg's method which involves maximal 
func\-tions. To make use of the corresponding results in \cite{stcon} we 
switch to the system
$$\psi _k^\alpha (x) = (2k!/\Gamma (k+\alpha +1))^{1/2} e^{-x^2/2}L_k^\alpha 
(x^2), \quad k\in {\bf N}_0, $$
which is obviously orthonormal on $L^2({\bf R_+}, d\mu _\alpha ),\; 
d\mu _\alpha (x) = x^{2\alpha +1}dx,\; \alpha > -1.$ For the sake of 
simplicity we write the norm on  $L^p({\bf R_+}, d\mu _\alpha )$ as
$$ \| F\| _p = \left( \int _0^\infty |F(x)|^pd\mu _\alpha (x) \right)^{1/p} .$$
We adopt 
the notion of the twisted generalized convolution  
on $L^1({\bf R}_+,d\mu _\alpha )$ from \cite{stcon} 
$$ F\times G(x) = \int _0^\infty  \tau _xF(y)G(y)\, d\mu _\alpha (y),$$
where the twisted generalized translation operator $\tau _x$ is given by
$$\tau _x F(y) = \frac{\Gamma (\alpha +1)}{\pi ^{1/2}\Gamma (\alpha +1/2)}
\int _0^\pi F((x,y)_\theta){\cal J}_{\alpha -1/2}(xy\sin 
\theta) (\sin \theta )^{2\alpha } \, d\theta ,$$
${\cal J}_{\beta }(x) = \Gamma (\beta +1) J_\beta (x)/(x/2)^\beta $, 
$J_\beta $ denoting the Bessel function of order $\beta $, and 
$$ (x,y)_\theta = (x^2 +y^2 -2xy \cos \theta)^{1/2}.$$
With respect to the system $\{ \psi _k^\alpha \}$ 
this convolution has the following transform property: 
if $F\sim \sum c_k\psi _k^\alpha $ and 
$F\times G \sim \sum c_kd_k\psi _k^\alpha $, then 
$G(x)\sim \Gamma (\alpha +1) \sum d_kL_k^\alpha (x^2)e^{-x^2/2}$. 

\medskip \noindent
If we set $ f(y^2)=F(y),\; g_\sigma (y^2)=G_\sigma (y)$, we see that 
$$|I_\sigma f(x^2)| = |F\times G_\sigma (x)|.$$
Apart from Lemma 2.1 the proof of Theorem 1.1 will be based on the fact 
that for $\alpha \ge 0$ and suitable $F$ and $G$ 
$$ |F\times G| \le |F| * |G| \, ,$$
which follows at once from the definition of the generalized Euclidean 
$*$-convolution 
$$F*G(x)=\int _0^\infty \tau^E_x F(y)G(y)\, d\mu _\alpha (y) \, $$
with associated generalized Euclidean translation 
$$\tau ^E_x F(y) = \int _0^\pi F((x,y)_\theta)\, d\nu _\alpha (\theta),\quad
d\nu _\alpha (\theta)=\frac{\Gamma (\alpha +1)}{\pi ^{1/2}\Gamma (\alpha +1/2)}
(\sin\theta)^{2\alpha}d\theta\,. $$
Therefore we restrict ourselves to fractional integrals defined 
via the generalized Euclidean convolution.

\medskip \noindent
\thm{
Let $1<p<q<\infty ,\; \alpha > -1/2,\;  
K_\sigma (x) = x^{2(\sigma -\alpha -1)}$. Then 
$$\|F*K_\sigma \| _q \le C \| F \| _p\, ,\quad
\frac1q =\frac1p -\frac{\sigma }{\alpha +1} \, .$$
}
By the above it is clear that Theorem 1.1 for the case $a=b=0$ follows from 
Theorem 2.2.

\medskip \noindent
{\bf Proof.} Following Hedberg \cite{pro} we want to estimate
$ F*K_\sigma (x)$ pointwise by a suitable maximal function which in this 
setting turns out to be (see Stempak \cite[p. 138]{stcon})
$$ F^*(x)=\sup_{\varepsilon>0}\varepsilon^{-(2\alpha+2)}
\int_0^\varepsilon\tau^E_x(|F|)(y)d\mu _\alpha (y) $$
with the usual boundedness property
$||F^*||_r\le C||F||_r,\quad 1<r\le \infty  . $ Now there holds 

\medskip \noindent
$$|F*K_\sigma (x)|
\le C \left( \int _0^\delta +\int _\delta ^\infty \right) \tau _x^E(|F|)(y) 
y^{2(\sigma -\alpha -1)} d\mu _\alpha (y) =J_1+J_2 ,$$
where $\delta >0$ will be chosen later appropriately. Clearly,
$$ J_1= \sum_{k=0}^\infty\int_{2^{-k-1}\delta}^{2^{-k}\delta} \ldots 
\leq
C \sum_{k=0}^\infty(2^{-k}\delta )^{2\sigma}(2^{-k}\delta )^{-2\alpha-2}
\int_{2^{-k-1}\delta}^{2^{-k}\delta}\tau^E_x(|F|)(s)s^{2\alpha+1}ds
$$
$$
\leq
C \delta^{2\sigma} \sum_{k=0}^\infty2^{-k2\sigma }(2^{-k}\delta)^{-2\alpha-2}
\int_0^{\delta2^{-k}}\tau^E_x(|F|)(s)s^{2\alpha+1}ds\leq
C\delta^{2\sigma }F^*(x).
$$
On the other hand, by H\"older's inequality,
$$
J_2=\int_\delta^\infty\tau^E_x(|F|)(s)s^{2(\sigma-\alpha-1)}s^{2\alpha+1}ds
\leq||\tau^E_x(|F|) ||_p \left (\int_\delta^\infty\,s^{2(\sigma-\alpha-1)p'}
s^{2\alpha+1}ds\right)^{1/p'}
$$
$$
\leq C\delta^{2\sigma -(2\alpha+2)/p}||F||_p\,\,,
$$
since $\tau^E_x$ are contractions on $L^p({\bf R_+}, d\mu _\alpha )$.
Hence 
$$ |F*K_\sigma (x)|
\le C\left(\delta^{2\sigma }F^*(x)+\delta^{2\sigma-(2\alpha+2)/p}
||F||_p\,\right) .
$$
Choosing $\delta=\left(F^*(x)/||F||_p\right)^{-p/(2\alpha+2)}$ (where $
||F||_p \neq 0$ is assumed) we  obtain
$$|F*K_\sigma (x)|
\leq C F^*(x)^{1-\sigma p/(\alpha+1)}\| F\| _p ^{\sigma p/(\alpha +1)}
$$
and then $\|  F*K_\sigma \| _q\leq C||F||_p$
due to the inequality  for the maximal function with $r=q\, (1-\sigma 
p/(\alpha+1))=p >1$.
This finishes the proof of Theorem 2.2.

\section{Extension to power weights}

\bigskip \noindent
The proof of Theorem 1.1 in the general case follows along the lines in Section 
2 from

\bigskip \noindent
\thm{
Let $\alpha >-1/2,\; 0<\sigma <\alpha +1$ and $ a<(\alpha +1)/p',\; 
b<(\alpha +1)/q,$ with $ a+b\ge 0$. If $1<p\le q<\infty $, then
$$
\| K_\sigma*F(x)\, x^{-2b}\| _q  \le  C \| F(x)\, x^{2a}\| _p , \quad 
\frac1q =\frac1p -\frac{\sigma -a-b}{\alpha +1}. 
$$
} 

\medskip \noindent 
{\bf Proof.} An equivalent version of the inequality above is
$$
\left( \int _0^\infty |Sf(x)|^q d\mu _\alpha (x)\right) ^{1/q}
\leq C\left( \int _0^\infty |f(x)|^p d\mu _\alpha (x)\right) ^{1/p}
$$
where $f(x) = x^{2a}F(x),\; Sf(x)=\int _0^\infty K(x,y)f(y)d\mu _\alpha (y)$ and
$$
K(x,y)=x^{-2b} \left( \int _0^\pi (x,y)_\theta ^{2(\sigma -\alpha -1)} 
d\nu _\alpha (\theta) \right) y^{-2a}.
$$
We first consider the case $p=q$. Then $\sigma = a+b$ and, 
therefore, the kernel $K(x,y)$ is homogeneous of degree $-(2\alpha+2)$:
$
K(rx,ry)=r^{-(2\alpha+2)}K(x,y), \; r>0.
$
It now suffices to check that (cf. Section 2 of \cite{ST})
$$
\int_0^\infty K(1,y)y^{-(2\alpha+2)/p}d\mu _\alpha (y)<\infty. 
$$
We first note that the function $K(1,y)$
has at most an integrable singularity at $y=1$, since for 
$1/2 \le y \le 2 $ we have
$$
\int_0^\pi\frac{(\sin\theta)^{2\alpha}d\theta }{((1-y)^2+4y\sin ^2 
(\theta/2))^{\alpha +1 -\sigma }}
\le C (1+|1-y|^{\sigma -1} )\, .$$ 

\smallskip \noindent
To deal with the singularity at 0 we note that for $y<1/2$ we have
$(1,y)_\theta  \approx 1$ and therefore $\int_0^\pi(1,y)_\theta^{2(\sigma 
-\alpha -1)}d\nu _\alpha(\theta ) \approx 1$
as well. Hence $K(1,y)\leq Cy^{-2a}$ and thus
$$
\int_0^{1/2}K(1,y)y^{-(2\alpha+2)/p}d\mu _\alpha (y)<\infty ,\quad a<(\alpha 
+1)/p' .$$
For $y>2$ we have $(1,y)_\theta \approx y$, $\int_0^\pi(1,y)_\theta^{2(\sigma 
-\alpha -1)} d\nu _\alpha (\theta)\approx y^{2(\sigma -\alpha -1)}$ and hence
$$
\int_2^{\infty}K(1,y)y^{-(2\alpha+2)/p}d\mu _\alpha (y)<\infty ,\quad 
b < (\alpha +1)/q,$$
since $\sigma < (\alpha +1)/p.$

\medskip \noindent 
Now, consider the case $1<p<q<\infty $. We use still another
equivalent version of the inequality to be proved, namely
$$
\int_0^\infty\int_0^\infty K(x,y)f(y)g(x)\, d\mu _\alpha (y)\, 
d\mu _\alpha (x)\leq C||f||_p||g||_{q'}\,\,,
$$
assuming for simplicity that $f$ and $g$ are nonnegative. Writing ${\bf R_+}\times
{\bf R_+}=D_1\cup D_2\cup D_3$ where
$$
D_1=\{y/2\le x\le 2y\},\,\, D_2=\{2y<x \} ,\,\, D_3=\{x<y/2\}
$$
it suffices to check for $i = 1, 2, 3$ that
$$
I_i=\int \! \! \!  \int _{D_i}K(x,y)f(y)g(x)\, d\mu _\alpha (y)\, 
d\mu _\alpha (x) \leq C \| f\| _p\| g\| _{q'}\,\,.
$$
Consider $I_1$ first. Since $a+b\ge 0$, for $x,y \in D_1$ and $\theta \in 
(0,\pi)$
we have
$$
\left( (x,y)_\theta \right) ^{2(a+b)}\leq Cx^{2(a+b)}\leq Cx^{2b} y^{2a};
$$
so 
$$ I_1\leq \int \! \! \! \int _{D_1}\widetilde K(x,y)f(y)g(x)\, d\mu _\alpha (y)
\, d\mu _\alpha (x)
$$
with $\widetilde K(x,y)=\tau^E_xK_{\sigma -a-b}(y)$\, . 
Hence we are reduced to showing that
$$
\int_0^\infty K_{\sigma -a-b} *f(x)\, g(x)\,d\mu _\alpha (x)\leq
||f||_p||g||_{q'},$$
which is implied by Theorem 2.2.

\medskip \noindent
To estimate $I_2$ and $I_3$ we need the following lemma.

\medskip \noindent
\la{
Let $\alpha > -1$ and $V_\delta f(x)=x^{2(\delta -\alpha -1)}
\int _0^xf(y)y^{-2\delta}d\mu _\alpha (y)$, $\delta <(\alpha +1)/p'$. Then 
$$||V_\delta f||_p\le C||f||_p\, ,\quad 
|V_\delta f(x)|\le C\, x^{-(2\alpha +2)/p}||f||_p\,.
$$
}

\medskip \noindent 
{\bf Proof.}
We have that $V_\delta f(x)=\int_0^\infty L(x,y)f(y)d\mu _\alpha (y),$ where 
$L(x,y)= x^{2(\delta-\alpha-1)}y^{-2\delta}$ for $y<x$ and equals 0 otherwise. 
Clearly $L$ is homogeneous of degree $-(2\alpha+2)$. 
Therefore the norm inequality is implied by
$$
\int_0^\infty L(1,y)y^{-(2\alpha+2)/p}d\mu _\alpha (y)=
\int_0^1y^{-2\delta-(2\alpha+2)/p+2\alpha+1} \, dy <\infty  $$
which is finite since $\delta < (\alpha +1)/p'$. For the pointwise estimate 
we simply write
$$
|V_\delta f(x)|\leq x^{2(\delta- \alpha -1)}\int_0^x|f(y)|\, y^{-2\delta}d\mu 
_\alpha (y) $$
$$
\leq C x^{2(\delta-\alpha -1)}\left(\int_0^xy^{-2\delta p'}d\mu _\alpha (y)
\right)^{1/p'}||f||_p \leq C x^{-(2\alpha+2)/p}||f||_p
$$
which finishes the proof of the lemma.

\medskip \noindent
Estimating $I_2$, we note that for $x,y\in D_2$ and all $\theta ,\; 0<\theta < 
\pi, $ there holds $x<2(x-y)\le 2((x-y)^2+2xy (1-\cos \theta ))^{1/2} =
2(x,y)_\theta $; hence
$$
I_2\leq C \int_0^\infty x^{2(\sigma -\alpha -1-b)}g(x)\left(\int_0^xf(y)\, y^{-2a}
d\mu _\alpha (y)\right)d\mu _\alpha (x)
$$
$$
= C \int_0^\infty g(x)x^{2(\sigma -a-b)}V_af(x)\, d\mu _\alpha (x)
\leq C ||g||_{q'} ||x^{2(\sigma -a-b)}V_af||_q\,.
$$
It now suffices to estimate the latter $L^q$ norm by $||f||_p$. We have
$$
\int_0^\infty|x^{2(\sigma -a-b)}V_a f(x)|^q d\mu _\alpha (x)
=\int_0^\infty|V_a f(x)|^p |V_a f(x)|^{q-p}x^{2(\sigma -a-b)q}
d\mu _\alpha (x)\,.
$$
Since $a<(\alpha+1)/p'$, by the pointwise estimate
in the above lemma 
$$
|V_a f(x)| \leq C x^{-(2\alpha+2)/p}||f||_p \, ;
$$
hence 
$$
|V_a f(x)|^{q-p}x^{2(\sigma -a-b)q}\leq C||f||_p^{q-p}
$$
due to the identity from assumptions. Hence 
$$
\int_0^\infty|x^{2(\sigma -a-b)}V_a f(x)|^qd\mu _\alpha (x)
\leq C ||f||_p^{q-p} \int_0^\infty|V_a f(x)|^pd\mu _\alpha (x)
\leq C ||f||_p^q\,.
$$
by the norm inequality from Lemma 3.2. The estimate of $I_3$ is similar and
we omit it. The proof of the theorem is complete. 

\bigskip \noindent
After submitting this paper the authors became acquainted with the following 
result of Kanjin and Sato [The Hardy-Littlewood theorem on fractional 
integration for Laguerre series, Proc. Amer. Math. Soc., to appear]:
$$\| I_\sigma f \| _{L^q_{v(\alpha q/2)}} \le C \| f \| _{L^p_{v(\alpha p/2)}} 
,\; \quad 0<\frac1q =\frac1p -\sigma ,\quad p>1,\quad \alpha \ge 0,$$
whereas for $-1 < \alpha <0$ there occurs $(1+\frac{\alpha }{2 })^{-1} <p,\, 
q<-\frac{2}{\alpha }$ as additional restriction.
Combining the above Theorem 1.1 in the case $a=b=\alpha =0$ with Kanjin's
transplantation theorem (cf. \cite{ST}) one at once recovers the above stated 
result of Kanjin and Sato .

\end{document}